\documentstyle{amsppt}
\magnification=\magstep 1
\vsize=8.5 truein
\hsize=6.5 truein
\tolerance=9000
\parindent=9pt
\parskip=9pt
\NoBlackBoxes\loadeusb
\loadbold
\loadmsbm
\loadmsam
\UseAMSsymbols
\pageno=1

\centerline{\bf EXTREMAL PROBABILISTIC PROBLEMS AND HOTELLING'S
$T^2$ TEST}

\centerline{\bf UNDER SYMMETRY CONDITION}

\centerline{by Iosif Pinelis}

\bigskip

\qquad$\topfoldedtext\foldedwidth{5.5 truein}{We consider
Hotelling's $T^2$ statistic for an arbitrary $d$-dimensional
sample. If the sampling is not too deterministic or inhomogeneous,
then under zero
means hypothesis, $T^2$ tends to $\chi^2_d$ in distribution. We are
showing that
a test for the orthant symmetry condition introduced by Efron can
be constructed
which does not essentially differ from the one based on $\chi^2_d$ and
at the same time
is applicable not only for large random homogeneous samples but for
all multidimensional
samples without exceptions. The main assertions have the form of
inequalities, not that
of limit theorems; these inequalities are exact representing the
solutions to certain extremal problems.
Let us also mention an auxiliary result which itself may be of
interest: $\chi_d-(d-1)^{\frac{1}{2}}$
decreases in distribution in $d$ to its limit $N(0,
\frac{1}{2})$.}$

\flushpar{\bf 1. Introduction}

Efron (1969), Eaton and Efron (1970) discovered that the
Hotelling's statistic
$$T^2=\overline X C^{-1}\overline X^T\tag 1.1$$
possesses a strong conservativeness property; here and further
$\overline X$ is the sample mean, $C$ is the sample covariance
matrix.

To describe it, put
$$R^2=\frac{T^2}{1+T^2}.\tag 1.2$$
Under some natural conditions, including zero means, both
$nT^2$ and $nR^2$ tend in distribution to $\chi^2_d$ as $n\to\infty$,
where $n$ is the
volume of the sample, $d$ is the dimension of the sample space.

Efron (1969), for $d=1$, Eaton and Efron (1970), for all $d$,
proved that
$$\bold Ef(n^{\frac{1}{2}} R)\le \bold E f(\chi_d)\tag 1.3$$
\vfill
\eject

\flushpar for $f(u)=u^{2m}$, $m=1, 2, \dots$; the only requirement
they put
on the
sample was the so-called orthant symmetry condition (we should not
like to discuss
now extensions of definition (1.1) applicable to an arbitrary, not
necessarily normal sample
when $C^{-1}$ can not exist).

Our main result here in the paper is that (1.3) holds for every
convex function
$f$ belonging to the class $C^2_{conv}$ of all even functions
having convex second
derivative; besides, the inequality
$$\bold P(n^{\frac{1}{2}}R\ge x)< c\cdot\bold P(\chi_d\ge x), \quad
x\ge 0,\tag 1.4$$
is being extracted from (1.3), with the best possible constant
$c=2e^3/9$.

The factor $2e^3/9$ can be found in Eaton (1974). Earlier, it was
proved in Eaton (1970) that
$$\bold E f(S_n)\le \bold Ef(\chi_1),\tag 1.5$$
where $S_n=\epsilon_1 x_1+\cdots+ \epsilon_n x_n$, $x^2_1+\cdots+
x^2_n=1$, $\epsilon_1, \dots , \epsilon_n$
are independent identically distributed (i.i.d.) random variables
(r.v.'s) with $\bold P(\epsilon_i=1)=\bold
P(\epsilon_i=-1)=\frac{1}{2}$, $i=1, \dots , n$; $f$
belongs to a certain subclass of the class $F$ of all
differentiable even
functions $f$ such that $[f^\prime (t+\Delta)+f^\prime
(t-\Delta)]/t$
is non-decreasing in $t>0$ for each real $\Delta$. Then, based on
(1.5), Eaton (1974)
obtained the inequality
$$\bold P(S_n\ge x)\le\frac{2e^3}{9} \frac{\varphi (x)}{x}
e^{-\frac{9}{2x^2}}
\left(1-\frac{3}{x^2}\right)^{-4}, \quad x>3^{\frac{1}{2}},$$
with $\varphi (x)=(2\pi)^{-\frac{1}{2}} e^{-\frac{x^2}{2}}$, and
stated the conjecture that
$$\bold P(S_n\ge x)\le\frac{2e^3}{9} \frac{\varphi(x)}{x}, \quad
x>2^{\frac{1}{2}}.$$

In this paper, we are proving that (1.5) holds for all $f$
of the class $F$ (which actually coincides with $C^2_{conv}$) and
that
$$\bold P(S_n\ge x)<\frac{2e^3}{9}(1-\phi(x))<\frac{2e^3}{9}
\frac{\varphi(x)}{x}\tag 1.6$$
for all $x>0$, where $\phi(x)=\int^x_{-\infty}\varphi (u)du$.

Moreover, we give analogues of (1.5), (1.6) for definite quadratic
forms
in $\epsilon_1, \dots , \epsilon_n$ (instead of linear forms such
as $S_n$) to which
(1.3), (1.4) are simple corollaries.

In section 2 below, we give our variants of strict definitions and
also representations of
$T^2$ and $R^2$ for any multidimensional sample.

In section 3, we reproduce the orthant symmetry condition by Efron.

In section 4, we present probabilistic problems related to (1.5),
(1.6).

In section 5, based on the monotonicity of the likelihood ratio, it
is shown that

\flushpar $\chi_d-(d-1)^{\frac{1}{2}}$ decreases in distribution to
its limit
$N\left(0, \frac{1}{2}\right)$ when $d\uparrow \infty$.

In section 6, the results of the preceding sections are used to
obtain (1.3), (1.4); the size
of the corresponding confidence region is discussing.

Section 7, Appendix, contains some proofs.

\bigskip

\flushpar {\bf 2. Hotelling's $T^2$ statistic for an arbitrary
distributed multidimensional sample.}

Let $X_1, \dots , X_n$ be a sample such as in the title of the
section,
$X_i\in\Bbb R^d$; we identify $\Bbb R^d$ with $M_{1d}$; $M_{nd}$
denotes
the set of all $n\times d$ real matrices. Consider
$$\overline X=\sum^n_{i=1} X_i/n, \,\,\, S=\sum^n_{i=1} X^T_i
X_i/n, \,\,\,
C=S-\overline X^T \overline X,$$
i.e., the sample mean, the matrix of second sample moments
and the sample covariance matrix, resp.; $^T$ means transposition.

If e.g. the sample is normally distributed, then $C^{-1}$ exists
with probability (w.p.) $1$, and $T^2$ statistic
is defined as
$$(n-1) \overline X C^{-1} \overline X^T.\tag 2.1$$

In general, it may be defined as in Eaton and Efron (1970):
$$T^2=\text{cot}^2\Theta, \tag 2.2$$
where $\Theta$ is the angle between vector
$$\nu=(1, \dots , 1)^T\in M_{n1}\tag 2.3$$
and the linear hull $L(X)$ of the columns of the matrix
$$X=(X^T_1 \dots X^T_n)^T,$$
$i$th line of which is $X_i$ (it is convenient to omit the factor
$n-1$ in (2.1), which corresponds to (1.1)). It is also reasonable
to define another
statistic,
$$R^2=\cos^2\Theta, \tag 2.4$$
i.e.,
$$R^2=\frac{T^2}{1+T^2},\tag 2.5$$
assuming that $R^2=1\Leftrightarrow T^2=\infty$ \quad (to specify
(1.2)).

We suggest a somewhat different approach to the general-case
definition of $T^2$
statistic, closer to (1.1). It leads to the same notions as in
(2.2),
(2.4). Namely, put
$$T^2=\lim\limits_{\epsilon\downarrow 0} \overline X(C+\epsilon
I)^{-1}\overline X^T, \tag 2.6$$

$$R^2=\lim\limits_{\epsilon\downarrow 0}\overline X(S+\epsilon
I)^{-1}\overline{X}^T,\tag 2.7$$
where $I$ is the unit matrix; the limits here always exist,
finite or infinite. (Replacing here $I$ by any strictly positively
definite symmetric
matrix (p.d.m.), one would obtain the equivalent definitions since
the matrix function
$A\to A^{-1}$ is monotone on the set of all p.d.m.'s; see Marshall
and Olkin (1979)). Indeed, considering
$$T^2_\epsilon=\overline{X} C^{-1}_\epsilon \overline{X}^T, \,\,\,
R^2_\epsilon=\overline{X} S^{-1}_\epsilon \overline{X}^T, \,\,\,
C_\epsilon=C+\epsilon I, \,\,\, S_\epsilon=S+\epsilon I,$$
we have
$$\align T^2_\epsilon R^2_\epsilon&=\overline X C^{-1}_\epsilon
\overline X^T\overline X
S_\epsilon^{-1}\overline X^T\\
&=\overline X
C^{-1}_\epsilon(S_\epsilon-C_\epsilon)S_\epsilon^{-1}\overline X^T
=T^2_\epsilon-R^2_\epsilon, \endalign$$
which implies (2.5) when using definitions (2.6), (2.7).

Consider the matrix
$$P=X(X^TX)^-X^T\tag 2.8$$
of the orthoprojector from $M_{n1}$ onto $L(X)$; $A^-$ denotes
a $g$-inverse matrix for $A$ (Rao (1965)). Taking into account
equalities
$\overline X=\nu^TX/n$, $S=X^TX/n$ and ``quasispectral"
representation for $X$
(Rao (1965), Ch. 1, appendix 6.1), it can be seen that
$$R^2=\overline{X}S^-\overline{X}^T,\tag 2.9$$
which is equivalent to
$$nR^2=\nu^TP\nu,\tag 2.10$$
i.e., definition (2.7) coincides with (2.4), and, in view of
(2.5), definitions (2.2) and (2.6) are equivalent too.

In contrast to (2.9), the equality $T^2=\overline{X}
C^-\overline{X}^T$ is not
always true; the simplest reason is that the right-hand side of it
cannot be infinite.

\bigskip

\flushpar {\bf 3. Orthant symmetry}.

This condition was defined in Efron (1969), Eaton and Efron (1970)
as
$$(X_1, \dots , X_n) \overset{D}\to= (\epsilon_1 X_1, \dots ,
\epsilon_n X_n)\tag 3.1$$
or, equivalently, 
$$X \overset{D}\to= \Delta_\epsilon X,\tag 3.2$$
where $\overset{D}\to=$ means equality in distribution,
$\Delta_\epsilon=diag\{\epsilon_1, \dots , \epsilon_n\}$
is the diagonal matrix with $\epsilon_1, \dots , \epsilon_n$ on its
diagonal,
$\epsilon_1, \dots , \epsilon_n$ are i.i.d.r.v.'s independent of 
$X$
$$\bold P(\epsilon_i=\pm 1)=\frac{1}{2}, \,\,\, i=1, \dots , n.$$

It was mentioned in Eaton and Efron (1970) that under the orthant
symmetry condition,
$$nR^2\overset{D}\to= \epsilon^TP\epsilon,\tag 3.3$$
where $\epsilon=(\epsilon_1, \dots , \epsilon_n)^T$. It may be also
deduced from
(2.8), (2.10), (3.2) and the equalities
$\epsilon=\Delta_\epsilon\nu$,
$\Delta^2_\epsilon=I$.

\bigskip

\flushpar {\bf 4. Some extremal probability problems.}

Let $f$ be a locally bounded Borel even  function on $\Bbb R$;
$A\in M_{nn}$; $A\ge 0$, i.e., $A$ is a nonnegatively definite
matrix;
$\epsilon_0, \epsilon_1, \dots , \epsilon_n$, $\xi_0, \xi_1, \dots , \xi_n$ are 
independent r.v.'s; $\epsilon=(\epsilon_1, \dots , \epsilon_n)^T$
as above; $\xi=(\xi_1, \dots , \xi_n)^T$; $\epsilon_i$ and $\xi_i$
are symmetrically distributed, $\bold P(\epsilon_i=1)=\frac{1}{2}$,
$\bold E\xi^2_i=1$, 
$i=0, 1, \dots , n$; $\xi_0$ does not coincide in distribution with
$\epsilon_0$; $x=(x_1, \dots , x_n)^T\in M_{n1}$.

In section 1, there were mentioned the class $C^2_{conv}$ of all
real even
functions $f$ on $\Bbb R$ having finite and convex (and hence
continuous) second derivative
$f^{\prime\prime}$ and the class $F$ of all real differentiable
even functions $f$ on $\Bbb R$
such that $[f^\prime (t+\Delta)+f^\prime (t-\Delta)]/t$ is
non-decreasing in $t>0$ for each $\Delta\in\Bbb R$.

{\bf Theorem 4.1.} The following statements are mutually
equivalent.

\quad (i) \ \ $\bold Ef((\epsilon^TA\epsilon)^{\frac{1}{2}})
\le \bold E f((\xi^T A\xi)^{\frac{1}{2}}) \forall A, \xi$;

\quad (ii) \ $\bold Ef(\epsilon_1 x_1+\cdots+ \epsilon_n x_n)\le
\bold Ef(\xi_1 x_1+\cdots+\xi_n x_n)\forall x, \xi$;

\quad (iii) $f\in C^2_{conv}$;

\quad (iv) \ $f\in F$;

\quad (v) \ \ $g_{c, f}\in C^2_{conv} \forall c\ge 0$, where
$g_{c, f}(u)=f((u^2+c)^{\frac{1}{2}}))$;

\quad (vi) \ $\exists \xi_0$ bounded w.p. 1 $\forall a, \,\,\,
b\in\Bbb R$

$$\bold E f(a+b\epsilon_0)\le \bold E f(a+b \xi_0); \tag 4.1$$

\quad (vii) \ $\forall \xi_0 \forall a, b\in\Bbb R$ (4.1) holds;

\quad (viii) $\exists \xi_0$ bounded w.p.1

$$\bold E f(\eta_0 + b\epsilon_0)\le \bold E f(\eta_0+b\xi_0)\tag
4.2$$
for every random variable $\eta_0$ bounded w.p.1 and independent of
$\epsilon_0, \xi_0$.

(Boundedness of $\xi_i$ is required nowhere but in (vi), (viii). Of
course, some
expectations may be infinite.)

Proof. See Appendix.

Define the class $C^2_{\uparrow, conv}$ of all the functions
$f\in C^2_{conv}$ which are convex on $[0, \infty)$ or,
equivalently, are non-decreasing on
$[0, \infty)$ and, besides, are nonnegative; in other words, $f\in
C^2_{conv}$
and $f(0)\ge 0, f^{\prime\prime}(0)\ge 0$; the additional
requirement of being
nonnegative does not diminish the generality of following
corollaries
4.2 - 4.4, but we do need it in corollaries 4.5, 4.6 and in some
places further.

{\bf Corollary 4.2.} If $f\in C^2_{\uparrow, conv}$, then
$\forall A, \xi, x$
$$\align \bold E f&((\eta^T A\eta)^{\frac{1}{2}})\le \bold E
f((\xi^T A\xi)^{\frac{1}{2}}),\\
\bold Ef&(\eta_1 x_1+\dots+ \eta_n x_n)\le \bold Ef(\xi_1
x_1+\cdots+ \xi_n x_n),\endalign$$
where $\eta=(\eta_1, \dots , \eta_n)^T$, $\eta_i$ are independent
r.v.'s,
$\bold E \eta_i=0$, $|\eta_i|\le 1$ \qquad w.p.1.

Proof. See Appendix.

{\bf Corollary 4.3.} If $P$ is an orthoprojector matrix, then
$$\bold Ef((\epsilon^TP\epsilon)^{\frac{1}{2}})\le \bold E
f(\chi_r),$$
where $r=$ rank $P$, $\chi_r=(\chi^2_r)^{\frac{1}{2}}$, $\chi^2_r$
is a r.v. with
$\chi^2_r$ distribution, $f\in C^2_{conv}$. If, moreover, $f\in
C^2_{\uparrow, conv}$,
then one can substitute here $\eta$ for $\epsilon$.

Proof. It suffices to note that
$\xi^T P\xi \overset{D}\to= \chi^2_r$ if
$\xi_i\sim N(0, 1)$.

This corollary in the case $f(u)=u^{2m}$, $m=1, 2, \dots , $
was proved, in essential, in (Eaton and Efron (1970),
corollary 6.1).

{\bf Corollary 4.4.} If $\xi_1$ is $N(0, 1)$ r.v. and $x^2_1
+\cdots+x^2_n=1$, then
$\forall f\in C^2_{conv}$
$$\bold E f(\epsilon_1 x_1+\cdots+\epsilon_n x_n)\le \bold E
f(\xi_1).$$
If, moreover, $f\in C^2_{\uparrow, conv}$, then we can substitute
$\eta$ for $\epsilon$.

Proof. Put $P=xx^T$ in corollary 4.3.

This corollary was proved in Eaton (1970, 1974) for $f\in F$ (and
so, in view of the
equivalence (iii) $\Leftrightarrow$ (iv) of theorem 4.1, for 
$f\in C^2_{conv}$) but under the additional
restriction $\bold E|f(T_n)|^{1+\delta}\le M$,
$\delta >0$, $T_n=\epsilon^T x$ or $\eta^Tx$, resp.

{\bf Corollary 4.5.} Under the conditions of corollary 4.3,
$$\align\bold P(\eta^TP \eta\ge u^2)&\le Q_r(u)\tag 4.3\\
&<\frac{2e^3}{9}\bold P(-\chi_r\ge u), \,\,\, u\ge 0,\tag
4.4\endalign$$
where
$$Q_r(u)=inf\{\bold E f(\chi_r)/f(u): \,\,\, f\in C^2_{\uparrow,
conv}, \,\,\, f(u)>0\}.$$

Proof. See corollary 4.3 for (4.3) and (4.13) below for (4.4).

{\bf Corollary 4.6}. Under the conditions of corollary 4.4,
$$\bold P(\eta_1 x_1+\cdots+\eta_n x_n\ge u)\le Q_1(u)/2\tag 4.5$$
$$\qquad\qquad\qquad\qquad <\frac{2e^3}{9}\bold P(\xi_1\ge u),
\,\,\, u\ge 0.\tag 4.6$$

Proof. Put $P=xx^T$ in corollary 4.5.

A statement close to (4.5) was given in Eaton (1974); (4.6) is an
improvement of corollaries 1, 2 in Eaton (1974) and of the
conjecture following those corollaries therein.

Let us provide further information on $Q_r(u)$ which, in
particular, contains (4.4), (4.6).

{\bf Proposition 4.7.}

$$\align Q_r(u)&=\min[1, r/u^2, W_r (u)]\tag 4.7\\
&=\cases 1&\text{if}\quad 0\le u\le r^{\frac{1}{2}},\\
r/u^2&\text{if}\quad r^{\frac{1}{2}}\le u\le \mu_r,\\
W_r(u)&\text{if}\quad u\ge \mu_r,\endcases\tag 4.8\endalign$$
where
$$\align W_r(u)&=inf\{(u-t)^{-3}\bold E(\chi_r-t)^3_+: \,\,\,
t\in (0, u)\},\\\mu_r&=\bold E \chi^3_r/\bold E\chi^2_r;\tag
4.9\endalign$$
besides,
$$\mu_r\in ((r+1)^{\frac{1}{2}}, \,\,\, (r+2)^{\frac{1}{2}}).\tag
4.10$$

Proof. See Appendix.

{\bf Proposition 4.8.} Constant $2e^3/9=4.463 \dots$ is the best
possible in (4.4), (4.6) in the sense that for each $r$
$$Q_r(u)\sim \frac{2e^3}{9}\bold P(\chi_r\ge u), \,\,\,
u\to\infty;\tag 4.11$$
here and in what follows $a\sim b$ means $a/b\to 1$.

Proof. See Appendix.

Consider the ratio
$$\Lambda_r(u)=Q_r(u)/\bold P(\chi_r\ge u),$$
and define
$$q=q(u)=q_r(u)=\int^\infty_u s^{r-1} e^{-s^2/2} I\{s>0\}ds$$
so that
$$\bold P(\chi_r\ge u)=q(u)/q(0);$$
here and further, $I\{\Cal A\}=1$ when $\Cal A$ is true, $I\{\Cal
A\}=0$
otherwise.

{\bf Proposition 4.9.}

$$\Lambda_r(u)<(2e^3/9)+3[\Cal J(a_u)-\Cal J(3)], \,\,\, u\ge
\mu_r, \tag 4.12$$
$$\Lambda_r(u)<2e^3/9, \,\,\, u\ge 0,\tag 4.13$$

where

$$\align a_u&=3q(u)q^{\prime\prime}(u)/q^\prime (u)^2>0,\\
\Cal J(a)&=6a^{-4}(e^a-1-a-a^2/2-a^3/6);\tag 4.14\endalign$$
$$a_u \uparrow 3 (u\ge \mu_r, \,\,\, u\uparrow \infty);\tag 4.15$$
$$\Cal J(a) \quad\text{increases in}\quad a>0\tag 4.16$$
(actually in all $a$, with $\Cal J(0)=\frac{1}{4}$, but we need not
this improvement).

Proof. See Appendix.

{\bf Proposition 4.10.} Function $\Lambda_1(u)$ increases in $u\ge
\mu_1$.

Proof. See Appendix.

The most non-trivial of propositions 4.7 - 4.10 is of course
proposition 4.9. It shows that $\Lambda_r(u)<\tilde\Lambda_r(u)$,
where $\tilde\Lambda_r(u)\left[=\frac{2e^3}{9}+3(\Cal J(a_u)
-\Cal J(3))\right]\uparrow\frac{2e^3}{9}, \,\,\, u\ge\mu_r$.
Apparently, $\Lambda_r(u)$ itself increases in $u\ge \mu_r$ for
each $r$
(to $2e^3/9$, in view of proposition 4.8), but we can prove this
fact only for $r=1$ (proposition 4.10);
a scheme of proof like that of proposition 4.10 may do for all $r$
though it becomes too complicated in the general case.
Nevertheless,
proposition 4.9 given here may often be more exact and useful than
that hypothetic
qualitative result.

\bigskip

\flushpar {\bf 5. The monotonicity of a likelihood ratio and a
stochastic majorization.}

Consider the family $\chi_r-(r-1)^{\frac{1}{2}}$, where $r$ is any
real number in $[1, \infty)$, $\chi_r$ has the density
$$p_r(u)=C_r u^{r-1} e^{-u^2/2}I\{u>0\},\tag 5.1$$
$C_r$ depends only on $r$; one can see that $(r-1)^{\frac{1}{2}}$
is
the mode of $\chi_r$.

We shall show that this family has monotone likelihood ratio and
hence
is stochastically monotone. Let $(E, \le)$ be any partially ordered
set.

We say that a family $(\xi_r: \,\,\, r\in E)$ of r.v.'s having
densities
$(p_r: \,\,\, r\in E)$ has monotone likelihood ratio (MLR) if the
implication
$$r\le d, \,\,\, r\in E, \,\,\, d\in E, \,\,\,
-\infty<s<t<\infty$$
$$\Rightarrow p_r(t)p_d(s)\ge p_d(t)p_r(s)\tag 5.2$$
is true. In the case $p_r(t)>0 \forall r, t$, this definition just
means that
$p_d/p_r$ is non-increasing on $\Bbb R$ when $r\le d$.

We say that a family of r.v.'s $(\xi_r: \,\,\, r\in E)$ with the
tails
$F_r(t)=\bold P(\xi_r\ge t)$ has monotone tail ratio (MTR) if the
implication
$$r\le d, \,\,\, r\in E, \,\,\, d\in E, \,\,\, -\infty
<s<t<\infty$$
$$\Rightarrow F_r(t)F_d(s)\ge F_d(t)F_r(s)\tag 5.3$$
is true.

A family $(\xi_r: \,\,\, r\in E)$ of r.v.'s is called
stochastically
monotone (SM) if
$$r\le d, \,\,\, r\in E, \,\,\, d\in E, \,\,\, t\in \Bbb R
\Rightarrow
F_r(t)\ge F_d(t);\tag 5.4$$
this definition of the stochastical monotonicity or, in other
words, monotonicity
in distribution, is generally accepted (see, e.g., Marshall and
Olkin (1979)).

{\bf Proposition 5.1.} If $(\xi_r)$ has MLR, then it has MTR.

Proof. If $(\xi_r)$ has MLR, and $-\infty<s<t<\infty$, then
$$F_r(t)[F_d(s)-F_d(t)]=\iint\limits_{s\le u<t\le v}p_r(v)p_d(u)dv du\ge$$
$$\ge \iint\limits_{s\le u<t\le v}p_d(v)p_r(u)dv
du=F_d(t)[F_r(s)-F_r(t)],$$
and so $F_r(t)F_d(s)\ge F_d(t)F_r(s)(r\le d)$.

{\bf Proposition 5.2.} If $(\xi_r)$ has MTR, then it is SM.

Proof. In (5.3), tend $s$ to $-\infty$.

{\bf Theorem 5.3.} The family $\{\chi_r-(r-1)^{\frac{1}{2}}: \,\,\,
r\ge 1\}$ has MLR.

Proof. Take $d\ge r\ge 1$ and put $a=(r-1)^{\frac{1}{2}}, \,\,\, 
b=(d-1)^{\frac{1}{2}}$. Then
$$(\log[p_d(u)/p_r(u)])^\prime=(r-d)u^2/(a+b)(u+a)(u+b)\le 0$$
if $u>-a$, hence $p_r(t)p_d(s)\ge p_d(t)p_r(s)$ when
$-a\le s<t$; if $s<-a$, then $p_r(s)=0$, so this case is trivial.

{\bf Corollary 5.4.} This family has MTR.

Proof. See proposition 5.1.

{\bf Corollary 5.5.} This family, $\{\chi_r-(r-1)^{\frac{1}{2}}:
\,\,\, r\ge 1\}$,
is SM.

Proof. See proposition 5.2.

{\bf Remark 5.6.} Consider $\xi_r=\chi_r-(r-1)^{\frac{1}{2}},
\,\,\, 1\le r<\infty$;
let $\xi_\infty$ have $N\left(0, \frac{1}{2}\right)$ distribution.
Then
$\xi_r\rightarrow \xi_\infty$ in distribution when $r\to\infty$,
and one may supplement the family
in statements 5.3 - 5.5 by $\xi_\infty$.

{\bf Remark 5.7.} If required additionally that inequality (5.2)
should be
strict for some $s, \,\,\, r\le d, \,\,\, r\ne d$ on a
positive-measure
set of values $t$ such that $t>s$, then the conclusions of
propositions 5.1, 
5.2 can be also improved so that inequalities (5.3), (5.4) become
strict for those values
of $s, r, d$. In particular, statements 5.3 - 5.6 can be improved
so that
inequalities (5.3), (5.4) become strict whenever $1\le r<d\le
\infty$,
$-(d-1)^{\frac{1}{2}}\le s<t<\infty$, for
$\xi_r=\chi_r-(r-1)^{\frac{1}{2}}$,
$1\le r<\infty$, and $\xi_\infty$ having $N\left(0,
\frac{1}{2}\right)$
distribution. E.g., for all $u\in\Bbb R, \,\,\, d\ge 1$,
$$\bold P(\chi_d-(d-1)^{\frac{1}{2}}\ge u)>1-\phi(2^{\frac{1}{2}}u),\tag 5.5$$
where $\phi$ is $N(0, 1)$ distribution function; this is an
improvement of the central limit
theorem for $\chi_r$.

\bigskip

\flushpar {\bf 6. Inequalities for the distribution of $T^2$ under
the orthant symmetry condition.}

It is known that if $X_1, \dots , X_n, \dots$ are i.i.d. random
vectors with zero
means and a finite non-degenerate matrix of second moments, then
the distribution of
$T^2$ (and hence that of $R^2$), is close to that of $\chi_d^2/n$
when
$n\to\infty$. Similarly to Efron (1969), Eaton and Efron (1970), we
are showing that some 
kind of conservativeness holds under the orthant symmetry condition
only; this conservativeness is
stronger than that in Efron (1969), Eaton and Efron (1970).

We permanently suppose that the orthant symmetry condition takes
place.

{\bf Theorem 6.1.} If $f\in C^2_{\uparrow, conv}$, then
$$\bold E f(n^{\frac{1}{2}} R)\le \bold E f (\chi_d),$$
where $R=(R^2)^{\frac{1}{2}}$. (When $f(u)=u^{2m}$, $m=1, 2, \dots
, $ this 
inequality was proved in Eaton and Efron (1970)).

Proof. See (3.3) and corollary 4.3.

The class $C^2_{\uparrow, conv}$ is sufficiently large to make it
possible
to obtain the following, main result of this paper.

{\bf Theorem 6.2.} For all $u\ge 0$
$$\bold P(n^{\frac{1}{2}}R\ge u)\le Q_d(u)<\frac{2e^3}{9}\bold P(\chi_d\ge u).\tag 6.1$$

Proof. See (3.3) and corollary 4.5.

One can find further information about $Q_d(u)$ in propositions 4.7 - 4.10.

Consider now $\delta$-quantiles $\tilde x_d(\delta)$ and
$x_d(\delta)$,
$0<\delta<1$, for $n^{\frac{1}{2}}R$ and $\chi_d$, resp., i.e.,
$$\bold P(x_d\ge x_d(\delta))=\delta, \,\,\, \tilde x_d(\delta)=inf\{x\in \Bbb R: \bold P(n^{\frac{1}{2}}R\ge x)\le
\delta\}.$$
In particular, $\bold P(n^{\frac{1}{2}}R\ge \tilde x_d(\delta))\ge
\delta>\bold P(n^{\frac{1}{2}}R>\tilde x_d(\delta))$.

{\bf Theorem 6.3.} If $\delta\le 0.5$, then
$x_d(\delta)>(d-1)^{\frac{1}{2}}$ and
$$\align \tilde x_d(\delta)&< x_d(\delta/c)\tag 6.2\\
&< x_d(\delta)+[x_d(\delta)-(d-1)/x_d(\delta)]^{-1}\log c\tag 6.3\\
&< x_d(\delta)+[x_d(\delta)-(d-1)^{\frac{1}{2}}]^{-1}\log c\tag 6.4\\
&< x_d(\delta)+o(1)\tag 6.5\\
&< x_d(\delta)(1+d^{-\frac{1}{2}}\cdot o(1))\tag 6.6\\
&< x_d(\delta)(1+o(1)),\tag 6.7\endalign$$
where $o(1)\rightarrow 0$ uniformly in $d, n, X_1, \dots , X_n$
when $\delta \downarrow 0$; $c=2e^3/9$.

Proof. See Appendix.

Thus this theorem means that the size $\tilde x_d(\delta)$ of
Hotelling's
criterion under the orthant symmetry condition can exceed the limit
size
$x_d(\delta)$, if can, only by some negligible value.
It is just that conservativeness property we spoke about.

Moreover, this theorem implies that the greater the
dimension is, the even better the situation becomes (the same
tendency
takes place when $\delta$ decreases). This is illustrated
by the following numerical data for $\delta=0.05$, where $x_\delta$
stands for $x_d(\delta)$,
$z_\delta=x_\delta+[x_\delta-(d-1)/x_\delta]^{-1}\cdot \log c$.

$$\vbox{\def\nh{\noalign{\hrule}}\offinterlineskip
\halign{\strut\vrule# \tabskip=.5em plus .5em
&\hfil
#\hfill&\vrule
#&\hfill
#\hfill&\vrule
#&\hfill
#\hfill&\vrule
#&\hfill
#\hfill&\vrule
#&\hfill
#\hfill&\vrule
#&\hfill
#\hfill&\vrule
#&\hfill
#\hfill&\vrule
#&\hfill
#\hfill&\vrule
# &\hfill
#\hfill
\hfill&\tabskip=0pt\vrule#\cr\nh
&$d$& &1& &2& &5& &10& &20& &50& &$\infty$&\cr\nh
&$x_\delta$& &1.96& &2.45& &3.33& &4.28& &5.61& &8.22&
&$d^{\frac{1}{2}}+1.16$&\cr
&$x_{\delta/c}$& &2.54& &3.00& &3.85& &4.78& &6.10& &8.69&
&$d^{\frac{1}{2}}+1.61$&\cr
&$z_\delta$& &2.72& &3.18& &4.03& &4.97& &6.28& &8.88&
&$d^{\frac{1}{2}}+1.80$&\cr\nh}}$$

{\bf 7. Appendix.}

Proof of theorem 4.1. It consists in checking the following
implications:

(vii) $\Rightarrow $ (vi) $\Rightarrow $ (viii) $\Rightarrow $ (iii)
$\Rightarrow $ (iv) $\Rightarrow $ (vii),

(iii) $\Leftrightarrow$ (v) $\Rightarrow $ (i) $\Rightarrow $ (ii)
$\Rightarrow $ (vii).

Let us proceed.

\qquad (vii) $\Rightarrow $ (vi). It is trivial.

\qquad (vi) $\Rightarrow $ (viii). Inequality (4.2) is an
``average" of (4.1).

\qquad (viii) $\Rightarrow $ (iii). For the beginning, suppose that
$f$ has continuous
4th derivative $f^{(4)}$. Then Taylor's expansion gives
$$(\bold Ef(a+b\xi_0)-f(a)-f^{\prime\prime}(a)b^2/2)\cdot
24b^{-4}\rightarrow f^{(4)}(a)\cdot\bold E\xi^4_0$$
when $b\rightarrow 0$. Using the analogous expansion with $\epsilon_0$
instead of $\xi_0$, one can
see that $f^{(4)}(a)(\bold E\xi^4_0-1)\ge 0$. But $\bold
E\xi^4_0>1$
since, as it was fixed, $\xi_0$ does not coincide with $\epsilon_0$ in
distribution and $\bold E\xi^2_0=1$. Thus
$f^{(4)}(a)\ge 0$, and so $f^{\prime\prime}$ is convex.
To check now the general case, put $f_m(u)=\bold E f(u+\eta_0/m),
m=1, 2, \dots ,$ where
$\eta_0$ is a bounded r.v. with a sufficiently smooth density. By
what has been already
proved, $f^{(4)}_m(u)\ge 0 \,\,\, \forall u\epsilon\Bbb R$. Consider
the operator $\Delta^3f(u)=f(u+3)-3f(u+2)+3f(u+1)-f(u)$. Then
$\forall u\ge 0 \exists \Theta\in (0, 3)\Delta^3
f_m(u)=f_m^{\prime\prime\prime}(u+\Theta)\ge 
f_m^{\prime\prime\prime}(u) \ge f_m^{\prime\prime\prime}(0)=0$
since
$f^{(4)}_m\ge 0$ and therefore $f_m^{\prime\prime\prime}$ is
non-decreasing.
Note that $f_m(u)\rightarrow f(u)$ almost everywhere (a.e.). Hence,
$\Delta^3f_m(u)$ is bounded a.e., and so is
$f_m^{\prime\prime\prime}(u)$.
By Helly's theorem, $f_m^{\prime\prime\prime}\rightarrow h$ weakly
on each compact in
$\Bbb R$, for some subsequence of the values of $m\to\infty$ and
some non-decreasing finite function $h$. Integrating by parts 3
times or, more exactly,
using Fubini theorem, it is easy to see that
$$f\in C^2_{conv}\Leftrightarrow f(u)=a+bu^2/2+\int\limits_{t\ge 0}
(|u|-t)^3_+ \nu_f(dt)/6\tag 7.1$$
for some real $a, b$ and a $\sigma$-finite nonnegative measure
$\nu_f$, where $u_+=\max(u, 0)$; besides, if $f\in C^2_{conv}$,
then
$a=f(0), b=f^{\prime\prime}(0), \forall u\ge 0 \nu_f([0,
u])=f^{\prime\prime\prime}(u+0)$,
where $f^{\prime\prime\prime}(u+0)$ is the right derivative of
convex function
$f^{\prime\prime}(u)$; (7.1) is but a kind of Taylor's expansion.
Further,
$$\int_{t\ge 0}(|u|-t)^3_+ df_m^{\prime\prime\prime}(t)\rightarrow
\int_{t\ge 0}(|u|-t)^3_t dh(t)$$
since $f_m^{\prime\prime\prime}\rightarrow h$ weakly on compacts.
Hence,
$\exists a, b$
$$f_m(0)\rightarrow a, \,\,\,
f_m^{\prime\prime}(0)\rightarrow b,$$
$$f(u)=a+bu^2/2+\int_{t\ge 0}(|u|-t)^3_+dh(t);$$
thus, $f\in C^2_{conv}$.

(iii) $\Rightarrow $ (iv). Essentially, it was proved in (Eaton
(1974), Lemma 1); we
need only to substitute $f^{\prime\prime\prime}(u+0)$ for
$f^{\prime\prime\prime}(u)$
in Eaton (1974).

(iv) $\Rightarrow $ (vii). Put
$g(x)=f(a+bx^{\frac{1}{2}})+f(a-bx^{\frac{1}{2}}), \,\,\, x\ge 0$.
Then, $f\in F\Leftrightarrow g^\prime$ is non-decreasing in $x>0$.
Therefore, $g$ is convex on $[0, \infty)$, and
$$2\bold Ef(a+b\epsilon_0)=g(1)\le \bold Eg(\xi^2_0)=2\bold
Ef(a+b\xi_0),$$
by Jensen's inequality.

(v) $\Rightarrow $ (iii). Take $c=0$.

(iii) $\Rightarrow $ (v). In view of (7.1), it is sufficient to
prove that
$\forall t\ge 0\forall c\ge 0$ the function
$g(u)=(z^{\frac{1}{2}}-t)^3_+$,
$z=u^2+c$, belongs to $C^2_{conv}$. Calculations show:

\qquad 1) $z=t^2\Leftrightarrow u=u_+$ or $u=u_-$, where
$u_\pm=\pm(t^2-c)^{\frac{1}{2}}$;
here it is necessary that $t^2\ge c$; further,
$$g^{\prime\prime\prime}(u_\pm+0)-g^{\prime\prime\prime}(u_\pm-0)
=6t^{-3}(t^2-c)^{\frac{3}{2}}\ge 0$$
is $t>0$; $u_\pm=0$,
$g^{\prime\prime\prime}(0+0)-g^{\prime\prime\prime}(0-0)=12>0$ when
$t=0$
(and, hence, $c=0)$;

\qquad 2) $z<t^2\Rightarrow g^{(4)}(u)=0$;

\qquad 3) $z>t^2\Rightarrow
g^{(4)}(u)=9cz^{-\frac{7}{2}}[(z-5t^2)c+4t^2z]$ 
$\ge 9cz^{-\frac{5}{2}}\cdot \min(z-t^2, 4t^2)\ge 0$ since $0\le
c\le z$.
Thus $g^{\prime\prime\prime}(u+0)$ is non-decreasing in $u\in\Bbb
R$, $g\in C^2_{conv}$.

\qquad (v) $\Rightarrow$ (i). Note that
$\epsilon^TA\epsilon=(\alpha\epsilon_1+\beta)^2+c$, where
$c\ge 0, \alpha, \beta$ do not depend on $\epsilon_1$. Function
$$\align h(u)&=f(((\alpha u+\beta)^2+c)^{\frac{1}{2}})+f(((-\alpha
u+\beta)^2+c)^{\frac{1}{2}})\\
&=g_{c, f}(\alpha u+\beta)+g_{c, f}(-\alpha u+\beta)\endalign$$
belongs to $C^2_{conv}$ as $g_{c, f}(u)$ does. Thus
$$2\bold Ef((\epsilon^TA\epsilon)^{\frac{1}{2}})=\bold Eh(\epsilon_1)\le \bold
Eh(\xi_1)
=2\bold Ef((\tilde\epsilon^t A\tilde\epsilon)^{\frac{1}{2}}),$$
where $\tilde\epsilon=(\xi_1, \epsilon_2, \dots , \epsilon_n)^T$, because, as we
have already proved,
(v) $\Rightarrow$ (iii) $\Rightarrow$ (iv) $\Rightarrow$ (vii).
Successively
replacing the remaining $\epsilon_i$'s by $\xi_i$'s, we are coming to
(i).

\qquad (i) $\Rightarrow$ (ii). Put $A=xx^T$.

\qquad (ii) $\Rightarrow$  (vii). Take $n=2$ and note that
$$\bold Ef(a+b\xi_0)=\bold Ef(-a+b\xi_0)=\bold Ef(a\epsilon_1+b\xi_0).$$

The theorem is proved.

Proof of corollary 4.2. If $f(u)=bu^2, \,\,\, b\ge 0$, then the
inequalities
we are proving are evident. In view of (7.1), it suffices to
consider the
case $f(u)=(|u|)-t)^3_+$, $t\ge 0$. In this case, calculations show
that
(cf. with the proof of statement (iii)$\Rightarrow$  (v) of theorem
4.1)
$g_{c, f}$ is convex on $\Bbb R$. Therefore, the function $h$ in
the proof of the 
implication (v)$\Rightarrow$  (i) of theorem 4.1 is convex. It
remains to use the
following inequality due to Hunt (1955), see also Eaton (1974):
$$\bold E G(\eta_1, \dots , \eta_n)\le \bold E G(\epsilon_1, \dots ,
\epsilon_n),$$
if $G$ is a function convex in each its argument.

In the following proofs we need some auxiliary results. Put
$$\align \gamma=\gamma(u)=\gamma_r(u)&=\bold
E(\chi_r-u)^3_+/C_r=\int^\infty_u(s-u)^3p_r(s)ds/C_r\\
&=\int^\infty_u(s-u)^3 s^{r-1}e^{-s^2/2}I\{s>0\}ds,\endalign$$
where $p_r$ is defined by (5.1), and, as in section 4,
$$\align q=q(u)=q_r(u)&=\int^\infty_u p_r(s)ds/C_r\\
&=\int^\infty_u s^{r-1}e^{-s^2/2}I\{s>0\}ds.\endalign$$

{\bf Lemma 7.1.} For all $r>0$, $j=0, 1, 2, 3, 4$
$$(-1)^j \gamma^{(j)}(u)\ge 0, \,\,\, u\in\Bbb R,\tag 7.2$$
$$(-1)^j\gamma^{(j)}(u)\sim 6u^{r-5+j}e^{-u^2/2}, \,\,\,
u\to\infty\tag 7.3$$
$$\gamma^{(3)}(u)=-6q(u),\tag 7.4$$
$$\gamma^{(4)}(u)=6u^{r-1}e^{-u^2/2}I\{u\ge 0\}=-6q^\prime(u),
\,\,\, u\ne 0,\tag 7.5$$
$$\align\gamma^{(5)}(u)=-6q^{\prime\prime}(u)&=6[u-(r-1)/u]q^\prime(u)=\\
&=-[u-(r-1)/u]\gamma^{(4)}(u), \,\,\, u>0,\tag 7.6\endalign$$
where $\gamma^{(j)}(u)=(d^j/du^j)\gamma(u)$, $j=1, 2, \dots ,
\gamma^{(0)}(u)=\gamma(u)$.

Proof. Equalities (7.4) - (7.6) are the results of simple
calculations.
By L'Hospital's rule,
$$\gamma(u)/(6u^{r-5}e^{-u^2/2})\sim \gamma^\prime
(u)/(6[(r-5)u^{r-6}-u^{r-4}]e^{-u^2/2})$$
$$\sim-\gamma^\prime(u)/(6u^{r-4}e^{-u^2/2})\sim\cdots\sim
\gamma^{(4)}(u)/(6u^{r-1}e^{-u^2/2})=1 (u\to\infty),$$
and so (7.3) is proved. In particular, $\gamma^{(j)}(u)\rightarrow
0 (u\to\infty)$.
This, together with (7.5), implies (7.2).

{\bf Lemma 7.2.} For each $r\ge 1$
$$q(u)q^{\prime\prime}(u)/q^\prime (u)^2\uparrow
1\quad\text{when}\quad u\ge (r-1)^{\frac{1}{2}}, \,\,\,
u\uparrow\infty,\tag 7.7$$
or, equivalently,
$$(u-(r-1)/u)^{-1}>-q(u)/q^\prime (u)\tag 7.8$$
$$>(u-(r-1)/u)/[(u-(r-1)/u)^2+1+(r-1)/u^2], \,\,\,
u>(r-1)^{\frac{1}{2}}.\tag 7.9$$

Proof. In view of (7.6), inequality (7.9) means that
$f(u)=q(u)q^{\prime\prime}(u)/q^\prime(u)^2$ increases in $u\ge
(r-1)^{\frac{1}{2}}$,
and (7.8) means that $f(u)<1$. On account of (7.3), (7.4), we see
that
$f(u) \to 1 (u\to\infty)$, and so

\qquad (7.7) $\Leftrightarrow$ (7.8) \& (7.9) $\Leftrightarrow$  (7.9).

\flushpar Write (7.9) as $g(u)>0$, $u>(r-1)^{\frac{1}{2}}$, where
we put
$$g(u)=q(u)+q^\prime (u) (u-(r-1)/u)/[(u-(r-1)/u)^2+1+(r-1)/u^2].$$
Taking into account (7.3), (7.6), we see that $g(u) \to 0$ as
$u\to\infty$ and
$$\align g^\prime (u)&=q^\prime
(u)[(1+(r-1)/u^2)^2+(u-(r-1)/u)(r-1)/u^3]\\
&\times [(u-(r-1)/u)^2+1+(r-1)/u^2]^{-2}<0, \,\,\,
u>(r-1)^{\frac{1}{2}}.\endalign$$
Hence $g(u)>0 \,\,\, \forall u>(r-1)^{\frac{1}{2}}$, and the proof
is completed.

Proof of proposition 4.7. Using (7.1) and arguments similar to
(2.8),
(2.9) of Eaton (1974), one can obtain (4.7).

Integration by parts gives, in view of (5.1),
$$\bold E\chi^j_r=(r+j-2)\bold E\chi^{j-2}_r, \,\,\, j>2-r.\tag
7.10$$
Thus, by Schwartz inequality,
$$\align \mu_r&=\bold E\chi^3_r/\bold E \chi^2_r>(\bold E
\chi^3_r/\bold E\chi_r)^{\frac{1}{2}}=(r+1)^{\frac{1}{2}},\\
\mu_r&< (\bold E \chi^4_r/\bold
E\chi^2_r)^{\frac{1}{2}}=(r+2)^{\frac{1}{2}}, \endalign$$
so (4.10) is proved.

Consider the function
$$\mu(t)=t-3\gamma (t)/\gamma^\prime(t).$$
We have
$$\align \mu^\prime
(t)&=(3\gamma(t)\gamma^{\prime\prime}(t)-2\gamma^\prime(t)^2)/\gamma^\prime (t)^2\\
&=2(\beta_3 \beta_1-\beta^2_2)/\beta^2_2>0,\endalign$$
by Schwartz inequality, where $\beta_i=\bold E(\chi_r-t)^i_+$. In
view
of (7.3), $\mu(t)\to\infty$ when $t\to\infty$, and so
$$t\leftrightarrow u=\mu(t)$$
is one-to-one increasing correspondence, under which, in
particular,
numbers $t\ge 0$ correspond to $u\ge \mu(0)=\mu_r$ (see (4.9)), and
vice versa. Put
$$F(t, u)=\bold E(\chi_r-t)^3_+/(u-t)^3=C_r\gamma (t)/(u-t)^3.$$
Then
$$\align (\partial/\partial t)F(t,
u)&=C_r(u-t)^{-4}(3\gamma(t)+(u-t)\gamma^\prime(t))\\
&=C_r(u-t)^{-4}\gamma^\prime (t)(u-\mu(t)), \,\,\, t<u,\endalign$$
and, taking into account that $\gamma^\prime (t)<0$, we deduce
$$\align W_r(u)&=F(\mu^{-1}(u), u)\\&=\min\{F(t, u): \,\,\, t\in
(-\infty, u)\}, \,\,\, u\ge \mu_r.\tag 7.11\endalign$$
In particular, $\forall u\ge \mu_r$
$$\align W_r(u)&\le F(0, u)=u^{-3}\bold E\chi^3_r \le u^{-2}\bold
E\chi^3_r/\mu_r=r/u^2\\
&\le r/\mu^2_r<1,\endalign$$
in view of (4.10), already proved, which implies (4.8) when $u\ge
\mu_r$.
If, conversely, $u\le \mu_r$, then
$$W_r(u)=F(0, u) \ge u^{-2}\bold E\chi^3_r/\mu_r=r/u^2,$$
and so, by (4.7), $Q_r(u)=\min [1, r/u^2]$. Thus (4.8) is
completely proved,  and so is the proposition.

Proof of proposition 4.8. By proposition 4.7,
$Q_r(u)=W_r(u)$ for $u\ge \mu_r$. And so (see (7.11))
$$Q_r(u)=F(\mu^{-1}(u), u)=-C_r\gamma^\prime
(\mu^{-1}(u))^3/(27\gamma(\mu^{-1}(u))^2), \,\,\, u\ge \mu_r.$$
Using (7.3), we see that for each $r$
$$Q_r(u)/\bold P(\chi_r\ge u)\sim
(2/9)(\mu^{-1}(u)^{r-2}/u^{r-2})\exp
\left\{\frac{u^2}{2}-\frac{\mu^{-1}(u)^2}{2}\right\}$$
$$\sim (2/9) \exp \{u(u-\mu^{-1}(u))\}\sim 2e^3/9, \,\,\,
u\to\infty.$$

Proof of proposition 4.9. First, note that assertion (4.15)
coincides with (7.7); (4.16) is a consequence of Taylor's
expansion.

Then, take for the beginning $u\ge \mu_r$. In view of (7.11),
$$Q_r(u)\le F(\tau, u), \tag 7.12$$
where
$$\tau=\tau(u)=u+3q/q^\prime,$$
and, as above, $q=q(u), \,\,\, q^\prime =q^\prime (u)$. Taylor's
expansion gives
$$\align &\gamma(\tau)=\gamma(u)-(u-\tau)\gamma^\prime
(u)+(u-\tau)^2\gamma^{\prime\prime}(u)/2\\
&-(u-\tau)^3 \gamma^{\prime\prime\prime}(u)/6+(u-\tau)^4 \int^1_0
\Theta^3\gamma^{(4)}(\tau+(u-\tau)\Theta)d\Theta/6,\\
e^a&=1+a+a^2/2+a^3/6+a^4\int^1_0 \Theta^3
e^{(1-\Theta)a}d\Theta/6.\endalign$$
Further, for $s<u$,

$\gamma^{(4)}(s)/\gamma^{(4)}(u)\le
[1-(u-s)/u]^{r-1}\exp\{(u^2-s^2)/2\}$

$<\exp\{(u-s)[-(r-1)/u+(u+s)/2]\}<\exp\{-(u-s)q^{\prime\prime}/q^
\prime\}$

$=\exp\{(1-\Theta)a_u\}, \,\,\, \Theta=(s-\tau)/(u-\tau), \,\,\,
a_u=3qq^{\prime\prime}/q^{\prime 2}$.

\flushpar Thus,
$$\align
\gamma(\tau)&<\gamma(u)-(u-\tau)\gamma^\prime(u)+(u-\tau)^2\gamma
^{\prime\prime}(u)/2\\
&-(u-\tau)^3
\gamma^{\prime\prime\prime}(u)/6+(u-\tau)^4\gamma^{(4)}(u)\Cal
J(a_u)/6,\endalign$$
$$\align &F(\tau, u)/\bold P(\chi_r\ge u)<\gamma\cdot
(u-\tau)^{-3}/q-\gamma^\prime\cdot (u-\tau)^{-2}/q\\
&+\gamma^{\prime\prime}\cdot
(u-\tau)^{-1}/(2q)-\gamma^{\prime\prime\prime}/(6q)+\gamma^{(4)}
\cdot (u-\tau)\Cal J(a_u)/(6q).\tag 7.13\endalign$$
Further,
$$\gamma^{(4)}\cdot (u-\tau)/(6q)=3,\tag 7.14$$
$$-\gamma^{\prime\prime\prime}/(6q)=1.\tag 7.15$$
Put $f_2=\gamma^{\prime\prime}+6q^2/q^\prime$. Then
$$f^\prime_2=\gamma^{\prime\prime\prime}+12q-6q^2q^{\prime\prime}
/q^{\prime 2}=6q(1-qq^{\prime\prime}/q^{\prime 2})>0$$
because of (7.7); hence, $f_2<0$, i.e.,
$$\gamma^{\prime\prime}\cdot (u-\tau)^{-1}/(2q)<1.\tag 7.16$$
Put $f_1=-\gamma^\prime-6q^3/q^{\prime 2}$. Then, by (7.7),
$$f^\prime_1=-\gamma^{\prime\prime}-18q^2/q^\prime+12q^3q^{\prime
\prime}/q^{\prime 3}>-\gamma^{\prime\prime}-6q^2/q^\prime=-f_2>0;$$
hence, $f_1<0$, i.e.,
$$-(u-\tau)^{-2}\gamma^\prime/q<2/3.\tag 7.17$$
Put $f_0=\gamma+6q^4/q^{\prime 3}$. Then, by (7.7),
$$f^\prime_0=\gamma^\prime+24q^3/q^{\prime 3}-18q^4q/q^{\prime
4}>\gamma^\prime+6q^3/q^{\prime 2}=-f_1>0;$$
hence, $f_0<0$, i.e.,
$$(u-\tau)^{-3}\gamma/q<2/9.\tag 7.18$$
Getting now (7.12) - (7.18), (4.15), (4.16) together, we obtain
$$\align \Lambda_r(u)&<2/9+2/3+1+1+3\Cal J(a_u)=\\
&=2e^3/9+3[\Cal J(a_u)-\Cal J(3)]<2e^3/9, \,\,\, u\ge
\mu_r.\endalign$$
Thus (4.12) and (4.13) for $u\ge \mu_r$ are proved.

Now consider the case $r^{\frac{1}{2}}\le u\le \mu_r$. Then (see
(4.7))
$$\Lambda_r(u)=rC^{-1}_r/(u^2q(u)), \,\,\,
(u^2q(u))^\prime=ug(u),$$
where we put $g(u)=2q(u)-u^re^{-u^2/2}$. We have
$$g^\prime(u)=[u^2-(r+2)]u^{r-1}e^{-u^2/2}<0$$
when $r^{\frac{1}{2}}\le u\le \mu_r(<(r+2)^{\frac{1}{2}}$; see
(4.10)). Thus,
$$g(r^{\frac{1}{2}})\le 0\Rightarrow \Lambda_r(u)\le
\Lambda_r(\mu_r)<\frac{2e^3}{9}, \,\,\, r^{\frac{1}{2}}\le u\le \mu_r,$$
in view of (4.12). If, on the contrary, $g(r^{\frac{1}{2}})>0$,
then
$$q_r(r^{\frac{1}{2}})>(r/e)^{r/2}/2.\tag 7.19$$
If $u_*\epsilon (r^{\frac{1}{2}}, \mu_r)$ is a root of the equation
$g(u)=0$, then
$$u^2_*q(u_*)=h(u_*)/2,$$
where we set $h(u)=u^{r+2}e^{-u^2/2}$; but
$$h^\prime(u)=-[u^2-(r+2)]h(u)/u>0$$
when $0\le u\le \mu_r(<(r+2)^{\frac{1}{2}})$; thus
$u^2_*q(u_*)>rq(r^{\frac{1}{2}})$;
hence
$$\max\{\Lambda_r(u): \,\,\, r^{\frac{1}{2}}\le u\le
\mu_r\}=\max\{\Lambda_r(r^{\frac{1}{2}}), \Lambda_r(\mu_r)\}.$$
By (4.8), (4.12), $\Lambda_r(\mu_r)<2e^3/9$. Let us prove that
$\Lambda_r(r^{\frac{1}{2}})<2e^3/9$ or, equivalently,
$q_r(0)<(2e^3/9)q_r(r^{\frac{1}{2}})$ under (7.19). It suffices to
show that
$$q_r(0)<\frac{e^3}{9}(r/e)^{\frac{r}{2}}.\tag 7.20$$
We shall do it by induction. This is easy when $r=1, 2$. If (7.20)
is true for some $r\ge 1$, then
$$\align &q_{r+2}(0)=rq_r(0)<r\cdot (e^3/9)(r/e)^{\frac{r}{2}}=\\
&=(e^3/9)[r/(r+2)]^{(r+2)/2}\cdot e\cdot
[(r+2)/e]^{(r+2)/2}<(e^3/9)
[(r+2)/e]^{(r+2)/2}.\endalign$$
Thus, (4.13) is proved in the case $r^{\frac{1}{2}}\le u\le \mu_r$.
Consider, finally, $0\le u\le r^{\frac{1}{2}}$. Then
$$\Lambda_r(u)=1/\bold P(\chi_r\ge u)\le
\Lambda_r(r^{\frac{1}{2}})<2e^3/9$$
as we have just shown. So (4.13) and the whole proposition are
completely proved.

Proof of proposition 4.10. Consider
$$\zeta (u)=Q^\prime_1(u)/(C_1 q^\prime_1(u)), \,\,\, u\ge \mu_1,$$
where $Q_1(u)=W_1(u)=F(\mu^{-1}(u), u)$ (see (4.8), (7.11)).
Let $t$ stand for $\mu^{-1}(u)$. Calculations show that
$$Q^\prime_1(u)=-\gamma^{\prime 3}/(27\gamma^2),$$
$$(d/dt)\log \zeta(u)=-3F^2/(\gamma\gamma^{\prime 3})\ge 0,$$
in view of (7.2), where
$F=\gamma^{\prime^2}-\gamma\gamma^{\prime\prime}, \,\,\,
\gamma=\gamma_1(t)$,
$\gamma^\prime=\gamma^\prime_1(t)$,
$\gamma^{\prime\prime}=\gamma^{\prime\prime}(t)$; moreover,
$$((((F/\gamma^{\prime\prime})^\prime
\gamma^{\prime\prime^2}/(\gamma^\prime
\gamma^{\prime\prime\prime}))^\prime 
\gamma^{\prime\prime\prime^2}/(\gamma^{\prime\prime}\gamma^{(4)})
)^\prime \gamma^{(4)^2}/(\gamma^{(5)}
\gamma^{\prime\prime\prime}))^\prime=6t^{-2}e^{-t^2/2}>0$$
when $t>0$, which, together with (7.2), (7.3), (7.5), leads to the
inequality $F>0$ and then to
$$(d/dt)\log\zeta(u)>0, \,\,\, t>0.$$
Now, in view of proposition 5.1 and remark 5.7, we obtain the
statement of
proposition 4.10.

Proof of theorem 6.3. Inequality (6.3) is a consequence of (6.1).
By (5.5), $x_d(\delta)>(d-1)^{\frac{1}{2}}$ when
$\delta\le 0.5$. Put $f(u)=cq(u+h)-q(u)$, where
$$c=2e^3/9, h=h(u)=(\log c)/[u-(d-1)/u], u>(d-1)^{\frac{1}{2}}.$$
Then
$$\align f^\prime(u)=-q^\prime (u)&+cq^\prime (u+h)\cdot
(1+h^\prime)>-q^\prime(u)+c q^\prime (u+h),\\
f^\prime (u) u^{1-d}e^{u^2/2}&>1-c\cdot
(1+h/u)^{d-1}\exp\left\{\frac{u^2}{2}-\frac{(u+h)^2}{2}\right\}\\
&>1-c\cdot \exp\{-[u-(d-1)/u]\cdot h\}=0.\endalign$$
Hence $f(u)<0 \,\,\, \forall u>(d-1)^{\frac{1}{2}}$. In view of
(6.1), this implies
(6.3), (6.4).

Using again (5.5), we see that if $\delta\downarrow 0$, then
$x_d(\delta)-(d-1)^{\frac{1}{2}}\rightarrow\infty$; hence we obtain
(6.5), (6.6).
Finally, (6.7) is trivial.
\vfill
\eject

\define\h{\hangindent=.4 truein\hangafter=1 }

\centerline{\bf REFERENCES}

\h Efron, Bradley. (1969). Student's $t$-test under symmetry
condition. {\it J. Amer.
Statist. Assoc. \bf 64} 1278-1302.

\h Eaton, Morris L. and Efron, Bradley. (1970). Hotelling's $T^2$
test under symmetry
conditions. {\it J. Amer. Statist. Assoc. \bf 65} 702-711.

\h Eaton, Morris L. (1974). A probability inequality for linear
combinations of bounded
random variables. {\it Ann. Statist. \bf 2} 609-614.

\h Eaton, Morris L. (1970). A note on symmetric Bernoulli random
variables.
{\it Ann. Math. Statist. \bf 41} 1223-1226.

\h Marshall, Albert W. and Olkin, Ingram. (1979). Inequalities:
theory of majorization and its applications. Academic Press.

\h Rao, C. Radhakrishna. (1965). Linear statistical inference and
its applications.
John Wiley \& Sons.

\h Hunt, G. A. (1955). An inequality in probability theory. {\it
Proc. Amer. Math. Soc. \bf 6} 506-510.

\bye